\documentclass{article}
\usepackage{amsfonts,amssymb,amsmath} 
\input{liemacs10.sty} 

\renewcommand{\phi}{\varphi}
\renewcommand{\mlabel}{\label}

\title{On the characterization of trace class representations and 
Schwartz operators} 
\author{Gerrit van Dijk,
  \begin{footnote}
{Math. Inst., Niels Bohrweg 1, 2333 CA Leiden, The Netherlands, 
dijk@math.leidenuniv.nl}   
  \end{footnote}
 Karl-Hermann Neeb, 
\begin{footnote}{
Department  Mathematik, FAU Erlangen-N\"urnberg, Cauerstrasse 11, 
91058-Erlangen, Germany; neeb@math.fau.de}
\end{footnote}\\
Hadi Salmasian,\begin{footnote}
{Department of Mathematics and Statistics,
University of Ottawa, 585 King Edward Ave., Ottawa, ON K1N 6N5,
Canada; hsalmasi@uottawa.ca}
\end{footnote}
 Christoph Zellner
\begin{footnote}{
Department  Mathematik, FAU Erlangen-N\"urnberg, Cauerstrasse 11, 
91058-Erlangen, Germany; zellner@math.fau.de}
\end{footnote}
} 
\begin{document} 

\maketitle

\begin{abstract} 
In this note we collect several characterizations 
of unitary representations $(\pi, \cH)$ of a finite dimensional Lie group $G$ 
which are trace class, i.e., for each compactly supported smooth 
function $f$ on $G$, the operator $\pi(f)$ is trace class. In particular 
we derive the new result that, for some $m \in \N$, 
all operators $\pi(f)$, $f \in C^m_c(G)$, are trace class. As a consequence 
the corresponding distribution character $\theta_\pi$ is of finite order. 
We further show $\pi$ is trace class if and only if 
every operator $A$, which is smoothing 
in the sense that $A\cH\subeq \cH^\infty$,  
is trace class and that this in turn is equivalent to the Fr\'echet space $\cH^\infty$ 
being nuclear, which in turn is equivalent to the realizability of the 
Gaussian measure of $\cH$ on the space $\cH^{-\infty}$ of distribution vectors. 
Finally we show that, even for infinite dimensional Fr\'echet--Lie groups, 
 $A$ and $A^*$ are smoothing if and only if $A$ is a Schwartz operator, 
i.e., all products of $A$ with operators from the derived representation 
are bounded. 
\end{abstract}

\section*{Introduction} 

Let $(\pi, \cH)$ be a (strongly continuous) unitary representation of the (possibly infinite dimensional) 
Lie group $G$ (with an exponential function).
Let $\cH^\infty$ be its subspace of smooth vectors. 
On this space we obtain by 
\[ \dd\pi(x)v = \frac{d}{dt}\Big|_{t=0} \pi(\exp tx)v \] 
the {\it derived representation} of $\g$ which we extend naturally 
to a representation of the enveloping algebra $U(\g)$, also denoted $\dd\pi$. 
We call an operator $A \in B(\cH)$ {\it smoothing} if $A\cH \subeq \cH^\infty$ 
(\cite{NSZ15}). A closely related concept is that of a {\it Schwartz operator}, 
which means that, for all $D_1, D_2 \in U(\g)$ (the enveloping algebra of the 
Lie algebra $\g$ of $G$), the sesquilinear form 
\[ (v,w) \mapsto \la A \dd\pi(D_2)v, \dd\pi(D_1) w \ra \] 
on $\cH^\infty$ extends continuously to $\cH \times \cH$ 
(\cite[Thm.~3.4, p.~349]{Ho77}, \cite{KKW15}). This note 
grew out of the question to understand the relation between smoothing 
and Schwartz operators. This is completely answered by Theorem~\ref{thm:2.4} which 
asserts, for any smooth representation of a Fr\'echet--Lie group $G$ 
and  $S \in B(\cH)$, the following are equivalent: 
\begin{itemize}
\item $S$ is Schwartz. 
\item $S$ and $S^*$ are smoothing. 
\item The map $G \times G \to B(\cH), (g,h) \mapsto \pi(g)S \pi(h)$ is smooth. 
\end{itemize}

Smoothing operators are of particular importance for 
unitary representations of finite dimensional Lie groups which are 
{\it trace class} in the sense that, for each $f \in C^\infty_c(G)$, 
the operator $\pi(f) = \int_G f(g)\pi(g)\, dg$ is trace class. 
Actually we show in Proposition~\ref{prop:1.6} that every 
smoothing operator is trace class if and only if $\pi$ is trace class. 
This connection was our motivation to compile various characterizations 
of trace class representations scattered in the literature, mostly without proofs  
(\cite{Ca76}). 
Surprisingly, this also led us to some new insights, 
such as the fact that, if $\pi$ is trace class, then 
there exists an $m \in \N$ such that 
all operators $\pi(f)$, $f \in C^m_c(G)$, are trace class. As a consequence,  
the corresponding distribution character $\theta_\pi$ is of finite order. 
This is contained in Theorem~\ref{thm:1.3} which collects various 
characterizations of trace class representations. 
One of them is that, for every basis $X_1, \ldots, X_n$ of $\g$ 
and $\Delta := \sum_{j = 1}^n X_j^2$, the positive selfadjoint operator 
$\1 - \oline{\dd\pi(\Delta)}$ has some negative power which is trace class. 
This is analogous to the Nelson--Stinespring 
characterization of CCR representations (all operators 
$\pi(f)$, $f \in L^1(G)$, are compact) by the compactness of the 
inverse of $\1 - \oline{\dd\pi(\Delta)}$. 
Locally compact groups for which all irreducible unitary representations are trace class 
have recently been studied in \cite{DD16}, and for a characterization 
of groups for which all irreducible unitary representations are CCR, we refer to 
\cite[Thm.~2]{Pu78}. 

In the measure theoretic approach to second quantization, the Fock space 
of a real Hilbert space is realized as the $L^2$-space for the Gaussian measure 
$\gamma$ on a suitable enlargement of $\cH$. Combining our 
characterization of trace class representations with 
results in \cite{JNO15}, we see that the trace class condition 
is equivalent to $\cH^\infty$ being nuclear, which in turn is equivalent 
to the realizability of the Gaussian measure on the dual space $\cH^{-\infty}$
of distribution vectors. 

{\bf Notation:} Throughout this article, $\N_0:=\N\cup\{0\}$. For a unitary representation $(\pi,\cH)$ of $G$, let $\oline{\dd\pi}(x)$ for $x\in\g $ denote the infinitesimal generator 
of the one-parameter group $t\mapsto\pi(\exp(tx))$ by Stone's Theorem.
Set
$\cD^n=\cD^n(\pi):=\bigcap_{x_1,\ldots,x_n\in\g}
\cD(\oline{\dd\pi}(x_1)\cdots\oline{\dd\pi}(x_n))
$ and
$\cD^\infty=\cD^\infty(\pi):=\bigcap_{n=1}^\infty\cD^n$.
\section{Characterizing trace class representations} 

In this section $G$ will be a finite dimensional Lie group and 
$\g$ will be the Lie algebra of $G$. We fix a basis $X_1, \ldots, X_n$ of $\g$ 
and consider the corresponding Nelson--Laplacian 
$\Delta := X_1^2 + \cdots + X_n^2$, considered as an element of 
the enveloping algebra $U(\g)$. We write $B_p(\cH)$ for the $p$th Schatten 
ideal in the algebra $B(\cH)$ of bounded operators on a Hilbert space~$\cH$ 
and $K(\cH)$ for the ideal of compact operators.

Recall that a unitary representation $(\pi,\cH)$ is called trace class if $\pi(f)\in B_1(\cH)$ for every $f\in C^\infty_c(G)$.
 For every  unitary representation 
$(\pi, \cH)$, the subspace $\cH^\infty$ of smooth vectors can naturally 
be endowed with a Fr\'echet space structure obtained from the 
embedding $\cH^\infty \to C^\infty(G,\cH), v \mapsto \pi^v$, 
where $\pi^v(g) = \pi(g)v$. Its range is the closed subspace of smooth equivariant 
maps in the Fr\'echet space $C^\infty(G,\cH)$.
This Fr\'echet topology on $\cH^\infty$ is identical to the topology obtained by the family of seminorms $\{\|\cdot\|_D\,:\,D\in U(\g)\}$, where
$\|v\|_D:=\|\dd\pi(D)v\|$ for $v\in\cH^\infty$.

\begin{lem} \mlabel{lem:1.1} If $(\pi, \cH)$ is a unitary representation 
of the Lie group $G$, then 
$\pi(f) \cH^\infty \subeq \cH^\infty$ for every $f \in C_c(G)$. 
\end{lem}

\begin{prf} In view of \cite[Thm.~4.4]{Ne10}, the representation 
$\pi^\infty$ of $G$ on the Fr\'echet space $\cH^\infty$ is smooth. 
Hence, for every $v \in \cH^\infty$ and $f \in C_c(G)$, 
the continuous compactly supported map 
\[ G \to \cH^\infty, \quad g \mapsto f(g)\pi(g) v\] 
has a weak integral $I$. Then, for every $w \in \cH$, 
\[ \la I,w \ra = \int_G f(g)\la \pi(g)v, w \ra\, dg = \la \pi(f)v, w \ra,\] 
and therefore $I = \pi(f)v \in \cH^\infty$. 
\end{prf}

\begin{lem} \mlabel{lem:1} Let $V$ be a Fr\'echet space,  $W$ be a metrizable vector space and 
$(\lambda_n)_{n \in \N}$ be a sequence of continuous linear maps $V \to W$ 
for which $\lambda(v)= \lim_{n \to \infty} \lambda_n(v)$ exists for every $v \in V$. 
Then $\lambda$ is continuous. 
\end{lem}

\begin{prf} Since $V$ is a Baire space and $W$ is metrizable, it follows from 
\cite[Ch. IX, \S 5, Ex. 22(a)]{Bou74} that 
the set of discontinuity points of $\lambda$ is of the first category, 
hence its complement is non-empty. This implies that $\lambda$ is continuous. 
\end{prf}

The following theorem generalizes \cite[Thm.~2.6]{Ca76} in a 
Bourbaki expos\'e of P.~Cartier which states the equivalence of (iii) and (v), 
but unfortunately without giving a proof or a reference to one. 

\begin{thm} \mlabel{thm:1.3}
  For a unitary representation $(\pi,\cH)$ of $G$, the following 
are equivalent: 
\begin{itemize}
\item[\rm(i)] There exists an $m \in \N$ such that $\pi(C^m_c(G)) \subeq B_1(\cH)$ and 
the corresponding map $\pi \:  C_c^m(G) \to B_1(\cH)$ is continuous. 
\item[\rm(ii)] $\pi(C^\infty_c(G)) \subeq B_1(\cH)$ and the map 
$\pi \: C^\infty_c(G) \to B_1(\cH)$ is continuous. 
\item[\rm(iii)] $\pi$ is a trace class representation, i.e., 
$\pi(C^\infty_c(G)) \subeq B_1(\cH)$. 
\item[\rm(iv)] $\pi(C^\infty_c(G)) \subeq B_2(\cH)$. 
\item[\rm(v)] There exists a $k \in \N$ such that 
$(\1 - \oline{\dd\pi(\Delta)})^{-k}$ is trace class. 
\end{itemize}
\end{thm}

\begin{prf} Let $D := \oline{\dd\pi(\Delta)}$, where 
$\dd \pi \:  U(\g) \to \End(\cH^\infty)$ denotes the derived representation,
extended to the enveloping algebra. Recall that $D$ is a non-positive 
selfadjoint operator on $\cH$ (\cite{NS59}). 

(i) $\Rarrow$ (ii) $\Rarrow$ (iii) $\Rarrow$ (iv) are trivial implications.

(iv) $\Rarrow$ (iii): 
According to the 
Dixmier--Malliavin Theorem (\cite[Thm.~3.1]{DM78}), 
we can write every $f \in C^\infty_c(G)$ as a finite sum of products 
$a* b$ with $a, b\in C^\infty_c(G)$. Hence the assertion follows from 
$B_2(\cH) B_2(\cH) \subeq B_1(\cH)$.

(iii) $\Rarrow$ (ii):
\begin{footnote}{The assertion in \cite[Prop.~1.4]{DD16} comes close to 
this statement but does not assert the continuity of the $B_1(\cH)$-valued 
map.}\end{footnote}
Let $(\delta_n)_{n \in \N}$ be a $\delta$-sequence in $C^\infty_c(G)$, 
i.e., $\int_G \delta_n(g)\, dg = 1$ and $\supp(\delta_n)$ converges to $\{\1\}$ 
in the sense that, for every $\1$-neighborhood $U$ in $G$, we eventually have 
$\supp(\delta_n) \subeq U$. Then 
$\delta_n * f \to f$ for every $f \in C^\infty_c(G)$ holds in $L^1(G)$ (and even in 
$C^\infty_c(G)$). 
For every $n \in \N$, the linear map 
\[\pi_n \:  C^\infty_c(G) \to B_1(\cH), \quad 
f \mapsto \pi(\delta_n) \pi(f) = \pi(\delta_n * f)\] 
is continuous because the linear maps 
\[ C^\infty_c(G) \to B(\cH),\quad  f \mapsto \pi(f) \quad \mbox{ and } \quad 
B(\cH) \to B_1(\cH), \quad A \mapsto \pi(\delta_n) A \] 
are continuous. Here we use that $\|\pi(\delta_n)A\|_1 \leq \|\pi(\delta_n)\|_1 \|A\|$. 

In view of Lemma~\ref{lem:1}, it suffices to show that, for every $f \in C^\infty_c(G)$, we 
have 
\[ \pi(f) = \lim_{n \to \infty} \pi_n(f) = \lim_{n \to \infty} \pi(\delta_n * f) \] 
holds in $B_1(\cH)$. Using the Dixmier--Malliavin Theorem (\cite[Thm.~3.1]{DM78}), 
we write $f = \sum_{j = 1}^k a_j * b_j$ with $a_j, b_j\in C^\infty_c(G)$. 
Then 
\[ \pi_n(f) = \pi(\delta_n * f) 
= \sum_{j = 1}^k \pi(\delta_n * a_j * b_j)
= \sum_{j = 1}^k \pi(\delta_n * a_j) \pi(b_j).\] 
Since the right multiplication maps  
$B(\cH) \to B_1(\cH), A \mapsto A \pi(b_j)$ 
are continuous and $\lim_{n \to \infty} \pi(\delta_n * a_j) = \pi(a_j)$ in $B(\cH)$, 
it follows that $\pi_n(f) \to \pi(f)$ for every $f \in C^\infty_c(G)$. 
Now the assertion follows from Lemma~\ref{lem:1}. 

(ii) $\Rarrow$ (v): Let $\Omega \subeq G$ be a compact $\1$-neighborhood 
in $G$ and 
\[ C^m_\Omega(G) := \{ f \in C^m(G) \: \supp(f) \subeq \Omega\} 
\quad \mbox{ for } \quad m \in \N_0 \cup \{\infty\}.\] 
Then $C^m_\Omega(G)$ is a Banach space for each $m\in \N_0$, and 
the Fr\'echet space $C^\infty_\Omega(G)$ is the projective limit of the 
Banach spaces $C^m_\Omega(G)$. Therefore the continuity 
of the seminorm $f \mapsto \|\pi(f)\|_1$ on $C^\infty_\Omega(G)$ implies the existence 
of some $m \in \N$ such that the map 
$\pi \: C^\infty_\Omega(G) \to B_1(\cH)$ extends continuously to 
$C^m_\Omega(G)$. This implies that $\pi(C^m_\Omega(G)) \subeq B_1(\cH)$. 

Next we observe that by an argument similar to the proof of a Lemma by M.~Duflo (\cite[Lemma~3.2.3, p.~250]{B72}), 
there exists for every $m\in \N$ a positive integer $k$, an 
open $\1$-neighborhood $U \subeq \Omega$ in $G$, and functions 
$\beta, \gamma \in C^m_c(U)$ such that
\begin{equation}
  \label{eq:duflo}
  (\1 - \Delta)^k \beta = \delta_\1 + \gamma, 
\end{equation}
where $\delta_\1$ is the Dirac distribution in~$\1$. Then
\[ \pi(\beta) 
= (\1 -D)^{-k} (\1 - D)^k \pi(\beta)
= (\1 -D)^{-k} \pi((\1 -\Delta)^k\beta)
= (\1 -D)^{-k} \big(\1 + \pi(\gamma)\big) \] 
holds as an identity of linear operators on $\cH^\infty$ 
(Lemma~\ref{lem:1.1}), and since both sides are bounded on $\cH$, we obtain 
\begin{equation}
  \label{eq:dag}
(\1 - D)^{-k} = \pi(\beta) - (\1 - D)^{-k} \pi(\gamma).
\end{equation}
By the preceding argument, both summands on the right are trace class, 
so that $(\1 -D)^{-k}$ is trace class as well. 

(v) $\Rarrow$ (i): For $f \in C^\infty_c(G)$, we have 
\begin{equation}
  \label{eq:rel1}
 \pi(f)
= (\1 - D)^{-k} (1-D)^k \pi(f) 
= (\1 - D)^{-k} \pi\big((\1-\Delta)^k f\big).
\end{equation}
Since the first factor on the right is trace class and $\pi((1-\Delta)^kf)\in B(\cH)$, it follows that 
$\pi(C^\infty_c(G)) \subeq B_1(\cH)$. Moreover, the continuity 
of the linear operator $(\1 -\Delta)^k \: C^{2k}_c(G) \to L^1(G)$ and 
the density of $C^\infty_c(G)$ in $C^{2k}_c(G)$ imply that the identity 
\eqref{eq:rel1} holds for all $f \in C^{2k}_c(G)$. We conclude that 
$\pi(C^{2k}_c(G)) \subeq B_1(\cH)$, and continuity of
the integrated representation 
$\pi \: L^1(G) \to B(\cH)$ implies 
that the corresponding map 
$C^{2k}_c(G) \to B_1(\cH)$ is continuous. 
\end{prf}

Along the same lines one obtains the following characterization 
of completely continuous representations (CCR) from \cite[Thm.~4.1]{NS59}. 

\begin{thm} {\rm(Nelson--Stinespring)} 
  For a  unitary representation $(\pi,\cH)$ of $G$, the following 
are equivalent: 
\begin{itemize}
\item[\rm(i)] $\pi(L^1(G)) \subeq K(\cH)$. 
\item[\rm(ii)] $\pi(C^\infty_c(G)) \subeq K(\cH)$. 
\item[\rm(iii)] $(\1 - \oline{\dd\pi(\Delta)})^{-1}$ is a compact operator. 
\end{itemize}
\end{thm}

\begin{prf} The equivalence of (i) and (ii) follows from the density 
of $C^\infty_c(G)$ in $L^1(G)$. We now use the same notation as in the preceding 
proof. 

(i) $\Rarrow$ (iii): From the relation 
\[ (\1 - D)^{-k} = \pi(\beta) - (\1 - D)^{-k} \pi(\gamma)\] 
we derive the existence of some $k \in \N$ for which $(\1 - D)^{-k}$ is compact, 
but this implies that $(\1 - D)^{-1}$ is compact as well.

(iii) $\Rarrow$ (ii): For $f \in C^\infty_c(G)$, we have 
\begin{equation}
  \label{eq:rel1b}
 \pi(f)
= (\1 - D)^{-1} (1-D) \pi(f) 
= (\1 - D)^{-1} \pi\big((\1-\Delta) f\big).
\end{equation}
Therefore the compactness of $(\1 - D)^{-1}$ implies (ii).
\end{prf}

\subsection*{Application to smoothing operators} 

\begin{defn} For a unitary representation $(\pi, \cH)$ of a Lie group 
$G$, an operator $A \in B(\cH)$ is called {\it smoothing} if 
$A \cH \subeq \cH^\infty$. We write $B(\cH)^\infty$ for the subspace of 
smoothing operators in $B(\cH)$. 
\end{defn}

It is shown in \cite[Thm.~2.11]{NSZ15} that for the class of Fr\'echet--Lie groups,
which contains in particular all finite dimensional ones, 
an operator $A$ is smoothing if and only if it is a smooth vector for the 
representation $\lambda(g)A := \pi(g)A$ of $G$ on $B(\cH)$. If 
$\pi$ is not norm continuous, then this representation is not continuous 
because the orbit map of the identity operator is not continuous, 
but it defines a continuous representation by isometries 
on the norm-closed subspace 
\[ B(\cH)_c := \{ A \in B(\cH) \: \lim_{g \to \1} \pi(g)A = A\}. \] 
By G\aa rding's Theorem, $\pi(f)\in B(\cH)^\infty$ for every $f\in C_c^\infty(G)$. Applying the Dixmier--Malliavin Theorem 
\cite[Thm.~3.3]{DM78} to the continuous representation 
$(\lambda, B(\cH)_c)$, 
we see that 
\begin{equation}
  \label{eq:eq1}
B(\cH)^\infty = \Spann \{ \pi(f)A \: f \in C^\infty_c(G), A \in B(\cH) \}.
\end{equation}
It follows in particular that 
all smoothing operators are trace class if $\pi$ is a 
trace class representation. Alternatively one can use the factorization 
\[ A = (\1 - D)^{-k} (\1 - D)^k A \] 
for every smoothing operator $A$ to see that $A$ is trace class because 
$(\1 - D)^{-k}$ is trace class for some~$k$. 

 From G\aa rding's Theorem we  
obtain another characterization of trace class representations: 

\begin{prop} \mlabel{prop:1.6} 
A unitary representation $(\pi, \cH)$ of $G$ is trace 
class if and only if $B(\cH)^\infty \subeq B_1(\cH)$, i.e., all smoothing 
operators are trace class. 
\end{prop}

\begin{prop} \mlabel{prop:1.7} If $(\pi, \cH)$ is a trace class representation of $G$, 
then the space of smoothing operators coincides with the subspace 
of smooth vectors of the unitary representation 
$(\lambda, B_2(\cH))$ defined by $\lambda(g)A := \pi(g)A$. 
\end{prop}

\begin{prf} Since the inclusion $B_2(\cH) \to B(\cH)$ is smooth, every 
$A \in B_2(\cH)^\infty$ has a smooth orbit map $G \to B(\cH), g \mapsto \pi(g)A$,
hence is smoothing. 

If, conversely, $A$ is smoothing, then \eqref{eq:eq1} shows that 
$A$ is a finite sum 
of operators of the form $\pi(f)B$, $f \in C^\infty_c(G)$, $B \in B(\cH)$.
Since $\pi \: C^\infty_c(G) \to B_2(\cH), f \mapsto \pi(f),$ is a continuous linear map by 
Theorem~\ref{thm:1.3}, the right multiplication 
map $B_2(\cH) \to B_2(\cH), C \mapsto CB$ 
is continuous, and the map 
$G \to C^\infty_c(G), g \mapsto \delta_g * f$ is smooth, the relation
 \[ \pi(g) \pi(f) B = \pi(\delta_g * f) B \] 
implies that $\pi(f)B$ has a smooth orbit map in $B_2(\cH)$. We conclude 
that the same holds for every smoothing operator.
\end{prf}

The equivalence of the statements in the first two parts of the following corollary can also be derived from  the vastly more general 
Theorem~\ref{thm:2.4}, but it may be instructive to see the direct argument 
for trace class representations as well. 

\begin{cor} 
\label{cor:1.8}
For a trace class representation 
of $G$ and $A \in B(\cH)$, the following are equivalent:
\begin{itemize}
\item[\rm (i)] $A$ is a Schwartz operator.
\item[\rm (ii)]
$A$ and $A^*$ are smoothing.
\item[\rm (iii)] $A\in B_2(\cH)$ and the map 
$\alpha^A \: G \times G \to B_2(\cH), (g,h) \mapsto \pi(g) A \pi(h^{-1})$ 
is smooth. 
\end{itemize}
\end{cor}

\begin{prf} (i) $\Rarrow$ (ii): 
If $A$ is Schwartz, then in particular the operators $A\dd\pi(D)$, $D\in U(\g)$, are bounded on $\cH^\infty$, and thus from \cite[Thm 2.11]{NSZ15} it follows that $A^*$ is smoothing. Furthermore, boundedness of $\dd\pi(D)A$ for every $D\in U(\g)$ entails in particular that $A\cH\subseteq \cD^\infty$, so that by \cite[Thm 2.11]{NSZ15} we obtain that $A$ is also smoothing. 

(ii) $\Rarrow$ (iii): Next assume that $A$ and $A^*$ are smoothing. Then 
Proposition~\ref{prop:1.6} implies that $A,A^*\in B_2(\cH)$, and Proposition~\ref{prop:1.7} 
implies that the maps 
\[ G \to B_2(\cH), \quad g \mapsto \pi(g)A \quad \mbox{ and } \quad 
  G \to B_2(\cH), \quad g \mapsto A\pi(g) \] 
are smooth. For the unitary representation of $G \times G$ on $B_2(\cH)$ 
defined by $\alpha(g,h)M := \pi(g)M \pi(h)^{-1}$ this implies that the matrix coefficient 
\[ (g,h) \mapsto  \la \alpha(g,h)A, A \ra 
= \la \pi(g)A\pi(h)^*, A \ra 
= \la \pi(g)A, A\pi(h)\ra  \] 
is smooth, so that $A$ is a smooth vector for $\alpha$ by 
\cite[Thm.~7.2]{Ne10}. 

(iii) $\Rarrow$ (i): 
Finally, assume that $A\in B_2(\cH)$ and the map $\alpha^A$ is smooth. Since the linear embedding $B_2(\cH)\to B(\cH)$ is continuous, the orbit map 
\begin{equation}
\label{orbGGBH}
G\times G\to B(\cH)\ ,\ (g,h)\mapsto \pi(g)A\pi(h)
\end{equation} is also smooth. 
From \cite[Lem 2.9]{NSZ15} and \cite[Lem 2.10]{NSZ15}, and by considering suitable partial derivatives 
at $(\1,\1)$
of the map \eqref{orbGGBH}, we obtain boundedness of
the operators 
\[ \dd\pi(X_1) \cdots \dd\pi(X_n) A \dd\pi(Y_1) \cdots \dd\pi(Y_m), \]
where $X_1, \ldots, X_n, Y_1, \ldots, Y_m \in \g$.
\end{prf}
For the last result of this section we need the following lemma, which appears in \cite[Thm 1.3(b)]{Ca76} without proof.
\begin{lem}
\label{lem:cartier1.3}
Let $(\pi,\cH)$ be a unitary representation of $G$ and let $\cH^{-\infty}$ denote the space of
distribution vectors, i.e., the anti-dual of $\cH^\infty$. Then every $\lambda\in \cH^{-\infty}$ is a sum of finitely many anti-linear functionals $\lambda_{D,v}\in\cH^{-\infty}$ of the form
$\lambda_{D,v}(w):=\langle v,\dd\pi(D)w\rangle$,
where $v\in \cH$ and $D\in U(\g)$.
\end{lem}
\begin{prf}
Continuity of  $\lambda_{D,v}$ is straightforward. Next fix $\lambda\in\cH^{-\infty}$. The map 
\[
\cH^\infty\to\cH^{U(\g)}\ ,\ 
v\mapsto(\dd\pi(D)v)_{D\in U(\g)}^{}
\]
is a topological embedding, where 
$\cH^{U(\g)}$ is equipped with the product topology. Thus by the Hahn--Banach Theorem, we can extend $\lambda$ to a continuous anti-linear functional on $\cH^{U(\g)}$. 
Since the continuous anti-dual of $\cH$ is identical to the continuous dual of the complex conjugate Hilbert space $\oline{\cH}$, and the continuous dual of a direct product is isomorphic to the direct sum of the continuous duals, we obtain that $\lambda=\sum_{i=1}^m
\lambda_{D_i,v_i}$ for some $D_i\in U(\g)$ and $v_i\in\cH$.  
\end{prf}
Let $\mathcal S(\pi,\cH)\subset B(\cH)$
denote the space of Schwartz operators of a unitary representation $(\pi,\cH)$. 
If $(\pi,\cH)$ is trace class, then from  Corollary \ref{cor:1.8} it follows that 
$\mathcal S(\pi,\cH)$ is the space of smooth vectors of
the unitary representation of $G\times G$ on the Hilbert space $B_2(\cH)$, defined by $\alpha(g,h)M:=\pi(g)M\pi(h)^{-1}$. In this case we equip  $\mathcal S(\pi,\cH)$ with the usual Fr\'echet topology of the space of smooth vectors. The next proposition characterizes the topological dual of $\mathcal S(\pi,\cH)$. 
\begin{prop}
\label{prp:linfunc}
Let $(\pi,\cH)$ be a trace class representation of $G$. 
Every continuous linear functional on the Fr\'echet space $\mathcal S(\pi,\cH)$ of Schwartz operators can be written as a sum of finitely many linear functionals 
\[ \lambda_{A,D}(T):=\tr(A\dd\pi(D)T\dd\pi(D')), \quad \mbox{  where } \quad 
A\in B_2(\cH), D,D'\in U(\g).\]   
\end{prop}

\begin{prf} From Corollary~\ref{cor:1.8} we know that 
the space $S(\pi,\cH)$ of Schwartz operators coincides with the space 
of smooth vectors of the unitary representation 
$(\alpha, B_2(\cH))$ given by $\alpha(g,h)A := \pi(g)A\pi(h)^{-1}$. 
For $x_1, \ldots, x_n,  y_1, \ldots y_m \in \g$ 
and every smooth vector $T$ for $\alpha$, we have 
\begin{align*}
& \dd\alpha((x_1,0) \cdots (x_n,0)(0,y_1) \cdots (0, y_m))T \\
&= (-1)^m \dd\pi(x_1) \cdots \dd\pi(x_n) T \dd\pi(y_m) \cdots \dd\pi(y_1).
\end{align*}
By Lemma \ref{lem:cartier1.3} if now follows that,
for every $T\in \mathcal S(\pi,\cH)$, we can write $\lambda$ as 
\begin{align*}
\lambda(T)&=
\sum_{i=1}^m\tr(\dd\pi(D_i)T\dd\pi(D'_i)A_i)
=
\sum_{i=1}^m\tr(A_i\dd\pi(D_i)T\dd\pi(D'_i)),
\end{align*}
where $A_i\in B_2(\cH)$ and $D_i\in U(\g )$ for
$1\leq i\leq m$.
\end{prf}

\subsection*{Nuclearity of the space of smooth vectors}

Combining Theorem~\ref{thm:1.3} with  \cite[Cor.~4.18]{JNO15} 
we obtain:  
\begin{prop} \mlabel{prop:3.18} For a 
unitary representation $(\pi, \cH)$, the following are equivalent: 
\begin{itemize}
\item[\rm(a)] $\pi$ is trace class. 
\item[\rm(b)] The Fr\'echet space $\cH^\infty$ is nuclear. 
\item[\rm(c)] There exists a measure $\gamma$ on the real dual space $\cH^{-\infty}$ 
of $\cH^\infty$, endowed with the $\sigma$-algebra generated by the evaluations 
in smooth vectors, whose Fourier transform is 
$\hat\gamma(v) = \int_{\cH^{-\infty}} e^{i \alpha(v)}\, d\gamma(\alpha) = e^{-\|v\|^2/2}$ 
for $v \in \cH^\infty$. 
\end{itemize}  
\end{prop}

The main idea in the proof of \cite[Cor.~4.18]{JNO15} is that 
$\cH^\infty$ coincides with the space of smooth vectors of 
the selfadjoint operator $\oline{\dd\pi(\Delta)}$ and that properties 
(b) and (c) can now be investigated in terms of the spectral resolution of this 
operator.  The equivalence of (a) and (b) is also stated in 
\cite[Thm.~2.6]{Ca76} without proof.

\section{Characterizing Schwartz operators} 

In this section we prove a characterization of Schwartz operators 
in terms of smoothing operators, namely that $S$ is Schwartz if and only if 
$S$ and $S^*$ are smoothing for any smooth unitary representation of a Fr\'echet--Lie 
group. 

We shall need the following result from interpolation theory 
(\cite[Prop.~9, p.~44]{RS75}: 

\begin{prop} \mlabel{prop:inter1} Let $\cH$ be a Hilbert space and 
$A, B$ be positive selfadjoint operators on $\cH$ with possibly unbounded 
inverses. Suppose that the 
bounded operator $T \in B(\cH)$ satisfies 
\[ T(\cD(A^2)) \subeq \cD(B^2) \quad \mbox{ with } \quad 
\|B^2 Tv\| \leq C \|A^2 v \| \quad \mbox{ for } \quad v \in \cD(A^2).\] 
Then $T(\cD(A)) \subeq \cD(B)$ with 
\[ \|B Tv\| \leq \sqrt{\|T\|C} \cdot \|A v \| \quad \mbox{ for } \quad v \in \cD(A).\] 
\end{prop}

We consider a smooth unitary representation $(\pi, \cH)$ of the 
(locally convex) Lie group $G$ and we assume that $G$ has a smooth exponential 
function. The next lemma provides an equivalent definition of Schwartz operators.



For $x \in \g$ and $n \in \N$, we consider the selfadjoint operator 
\[ N_{x,n} := \1 + (-1)^n\oline{\dd\pi}(x)^{2n} \geq \1.\] 
Note that \cite[Lemma~4.1(b)]{NZ13} implies that $N_{x,n}$ coincides 
with the closure of the operator $\1 + (-1)^n \dd\pi(x)^{2n}$ 
on $\cH^\infty$. 

\begin{prop} \mlabel{prop:1} If $S \in B(\cH)$ is a smoothing operator whose adjoint 
$S^*$ is smoothing as well, then $S$ is a Schwartz operator, i.e., 
for $D_1, D_2 \in U(\g)$, the operators $\dd\pi(D_1) S \dd\pi(D_2)$ 
defined on $\cH^\infty$ are bounded, i.e., extend to bounded operators on~$\cH$. 
\end{prop}

\begin{prf}  Since $U(\g)$ is spanned by the elements of the 
form $x^n$, $x \in \g$, we have to show that, 
for $x,y \in \g$ and $n,m \in \N_0$, 
the operator $\dd\pi(x)^n S \dd\pi(y)^m$ is bounded. 
From \cite[Thm.~2.11]{NSZ15} we know that 
$N_{x,n} S$ is bounded, and from
\cite[Lem.~2.8(a)]{NSZ15} it follows that
$S N_{y,m}$ is bounded on $\cD(N_{y,m})$. 
Next we observe that the operators 
\[ T := N_{x,n} S, \quad A := N_{y,m}^{-1/2} \quad \mbox{ and } \quad 
B := N_{x,n}^{-1/2} \] 
are all bounded. Writing $v \in \cH$  as 
$v = N_{y,m} w$, we obtain the estimate 
\[ \|B^2 T v\| 
= \| N_{x,n}^{-1} T N_{y,m} w \| 
= \| S N_{y,m} w \| 
\leq \| S N_{y,m}\| \|N_{y,m}^{-1}v \| 
= \| S N_{y,m}\| \|A^2v \|.\] 
Therefore Proposition~\ref{prop:inter1} implies that for $c:=\|N_{x,n}S\|^{1/2}\|SN_{y,m}\|^{1/2}$ we have 
\[ \|N_{x,n}^{1/2} Sv\| = \|B Tv\| \leq c \|Av\| = c \| N_{y,m}^{-1/2} v\| \quad 
\mbox{ for } \quad v \in \cD(A) = \cR(N_{y,m}^{1/2}) = \cH.\] 
For $ v= N_{y,m}^{1/2}w$, this leads to 
\[ \|N_{x,n}^{1/2} S N_{y,m}^{1/2} w\| \leq c \| w\|
\ \text{ for }w\in \cD(N_{y,m}^{1/2})
,\]
so that $N_{x,n}^{1/2} S N_{y,m}^{1/2}$ is bounded
on $\cD(N_{y,m}^{1/2})$. 
As $N_{x,n}^{-1/2}\dd\pi(x)^n$ is bounded, it follows that the following 
operator is bounded: 
\begin{align*}
&(N_{x,n}^{-1/2} \dd\pi(x)^n)^*(N_{x,n}^{1/2} S N_{y,m}^{1/2})(N_{y,m}^{-1/2} \dd\pi(y)^m) \\
&\supeq (\dd\pi(x)^n)^* (N_{x,n}^{-1/2} N_{x,n}^{1/2} S N_{y,m}^{1/2})(N_{y,m}^{-1/2} \dd\pi(y)^m) \\
&= (\dd\pi(x)^n)^* S  \dd\pi(y)^m \supeq (-1)^n \dd\pi(x)^n  S \dd\pi(y)^m, 
\end{align*}
and this implies the boundedness of $\dd\pi(x)^n S \dd\pi(y)^m$, more precisely
\begin{align}
\|\dd\pi(x)^n S \dd\pi(y)^m\| \leq \|N_{x,n}^{1/2} S N_{y,m}^{1/2}\| \leq \|N_{x,n}S\|^{1/2}\|SN_{y,m}\|^{1/2}. \label{prop1eq1}
\end{align}
\end{prf}

We now consider the representation of $G \times G$ on $B(\cH)$ by
\begin{align*}
 \alpha(g,h) A &:= \pi(g)A \pi(h)^{-1} = \lambda(g)\rho(h)A, \\ 
\lambda(g)A &= \pi(g)A, \qquad \rho(g)A = A \pi(g)^{-1}.
\end{align*}

\begin{rem} \mlabel{rem:1} (a) Suppose that $A$ is a continuous vector for 
the left multiplication representation $\lambda$ and also for the
right multiplication action $\rho$. Then 
\begin{align*}
 \|\pi(g)A\pi(h) - A \| 
&\leq \|\pi(g)A \pi(h) - A \pi(h)\| + \|A\pi(h) - A\|\\
&\leq \|\pi(g)A -A\| + \|A\pi(h) - A\|
\end{align*}
implies that $A$ is a continuous vector for $\alpha$. 

We write $B(\cH)_c(\alpha)$ for the closed subspace of $\alpha$-continuous vectors 
in $B(\cH)$ and note that since $\alpha$ acts by isometries, 
it defines a continuous action of $G$ on the Banach space 
$B(\cH)_c(\alpha)$. 

(b) Suppose that $A$ is a $C^1$-vector for $\lambda$ and $\rho$ 
and $x,y \in \g$. 
Since all operators $\pi(\exp tx)A, A \pi(\exp ty)$ are contained in 
$B(\cH)_c(\alpha)$, the closedness of $B(\cH)_c(\alpha)$ implies that 
$\oline{\dd\pi}(x)A$ and $A \oline{\dd\pi}(y)$ are also 
$\alpha$-continuous. 

We claim that $A$ is a $C^1$-vector for $\alpha$. In fact, the map 
\[ F \: G \times G \to B(\cH),\quad 
F(g,h) := \alpha(g,h)A = \pi(g)A \pi(h)^{-1} \] 
is partially $C^1$, so that its differential $\dd F \: TG \times TG \to B(\cH)$ 
exists. This map is given by 
\begin{align*}
\dd F(g.x, h.y) 
&= \pi(g) \oline{\dd\pi}(x)A \pi(h)^{-1} -\pi(g) A\oline{\dd\pi}(y) \pi(h)^{-1} \\
&= \pi(g)\big(\oline{\dd\pi}(x)A-A\oline{\dd\pi}(y)\big) \pi(h)^{-1}.
\end{align*}

Since $\alpha$ defines a continuous action on 
$B(\cH)_c(\alpha)$, the continuity of $\dd F$ 
follows from the continuity of the corresponding  linear map 
\[ \g \times \g \to B(\cH), \quad 
(x,y) \mapsto \oline{\dd\pi}(x)A-A\oline{\dd\pi}(y),\] 
which follows from the assumption that $A$ is a $C^1$-vector 
for $\lambda$ and $\rho$.
This shows that $\dd F$ is continuous and hence that $A$ is a $C^1$-vector for~$\alpha$. 
\end{rem}

\begin{thm} \mlabel{thm:2.4} For a smooth unitary representation of a Fr\'echet--Lie group 
$G$ and $S \in B(\cH)$, the following are equivalent: 
  \begin{itemize}
  \item[\rm(i)] $S$ and $S^*$ are smoothing. 
  \item[\rm(ii)] $S$ is Schwartz. 
  \item[\rm(iii)] $S$ is a smooth vector for $\alpha$, i.e., 
the map $G \times G \to B(\cH)$, \break $(g,h) \mapsto \pi(g)S\pi(h)^{-1}$ is smooth. 
  \end{itemize}
\end{thm}

\begin{prf} That (i) implies (ii) is Proposition~\ref{prop:1}. 

(ii) $\Rarrow$ (iii): For $D \in U(\g_\C)$, the operators 
$S\dd\pi(D)$ and $S^* \dd\pi(D)$ on $\cH^\infty$ are bounded, 
so that \cite[Thm.~2.11]{NSZ15} implies that $S$ and $S^*$ are smoothing, 
hence in particular $C^1$-vectors for $\alpha$ by 
Remark~\ref{rem:1} and 
\[ \oline{\dd\alpha}(x,y)S = \oline{\dd\pi}(x)S - S\oline{\dd\pi}(y).\] 
It follows in particular that $\oline{\dd\alpha}(x,y)S$ is Schwartz as well (because $S$ is Schwartz if and only if $\dd\pi(D_1)S\dd\pi(D_2)$ is bounded on $\cH^\infty$ for every $D_1,D_2\in\cH^\infty$). 
Thus we obtain inductively that $S \in \cD^n(\alpha)$ for every $n \in \N$. 
Since $G$ is Fr\'echet, \cite[Thm.~1.6(ii), Cor.~1.7]{NSZ15} now imply that 
$S$ is a smooth vector for $\alpha$. 

(iii) $\Rarrow$ (i) follows from the characterization of smoothing 
operators (\cite[Thm.~2.11]{NSZ15}). 
\end{prf}

If the Lie group $G$ is only assumed to be metrizable, the additional 
quantitative information from Proposition~\ref{prop:1} can still be used to obtain 
the equivalence of (i) and (iii) in the preceding theorem. This is done in Theorem \ref{thm:2.5} below. First we need a lemma.

\begin{rem}
\mlabel{rmk-hadi-xy}
Let $(\pi,\cH)$ be a smooth unitary representation of a 
(locally convex) Lie group $G$ 
with a smooth exponential map.
Let $S\in B(\cH)$ be a Schwartz operator, and set $A:=\dd\pi(D_1)S\dd\pi(D_2)$ with domain $\cH^\infty$, where $D_1,D_2\in U(\g)$. 
Then $A$ is bounded, and therefore $\oline{A}\in B(\cH)$. We now show that $\oline A(\cH)\subseteq\cD(\oline{\dd\pi}(x)^2)$. Indeed for 
$v\in \cH$, if $(v_n)_{n\in\N}\subseteq \cH^\infty$ is a sequence such that $\lim_{n\to\infty}v_n=v$ in $\cH$, then 
$\lim_{n\to \infty} Av_n=\oline Av$ and from boundedness of $\dd\pi(x)^2A$ with domain $\cH^\infty$ (recall that $S$ is Schwartz) 
it follows that the sequence $(\oline{\dd\pi}(x)^2Av_n)_{n\in\N}$ is convergent. But $\oline{\dd\pi}(x)^2$ is closed, hence $\oline Av\in\cD(\oline{\dd\pi}(x)^2)$.

\end{rem}

\begin{thm}
\label{thm:2.5}
Let $(\pi,\cH)$ be a smooth unitary representation of the Lie group $G$ and assume that $\g$ is metrizable. For $S \in B(\cH)$, the following are equivalent:
  \begin{itemize}
  \item[\rm(i)] $S$ and $S^*$ are smoothing. 
  \item[\rm(ii)] $S$ is a smooth vector for $\alpha$, i.e., 
the map $G \times G \to B(\cH)$, \break $(g,h) \mapsto \pi(g)S\pi(h)^{-1}$ is smooth. 
  \end{itemize}
\end{thm}
\begin{prf}
(ii) $\Rarrow$ (i) follows from \cite[Thm.~2.11]{NSZ15}. Now assume that $S$ and $S^*$ are smoothing. Then $S$ is a Schwartz operator by Proposition \ref{prop:1}. According to \eqref{prop1eq1} we have
\begin{align}
\|\dd\pi(x)^n S \dd\pi(y)^m\| &\leq \|N_{x,n}S\|^{1/2}\|SN_{y,m}\|^{1/2} \leq \textstyle{\frac{1}{2}}(\|N_{x,n}S\|+\|SN_{y,m}\|) \nonumber\\
& \leq \|S\| + \textstyle{\frac{1}{2}}(\|\dd\pi(x)^{2n}S\|+\|S\dd\pi(y)^{2m}\|).\label{metrthmeq1}
\end{align}
By \cite[Thm.~2.11]{NSZ15} the map $g\mapsto \lambda(g)S=\pi(g)S$ is smooth with 
$$\overline{\dd\lambda}(x_1)\cdots \overline{\dd\lambda}(x_n)S = \dd\pi(x_1)\cdots \dd\pi(x_n)S.$$
In particular, $\g^k \to B(\cH), (x_1,\dots,x_k) \mapsto \dd\pi(x_1)\cdots \dd\pi(x_k)S$ is $k$-linear and continuous. Similarly the smoothness of $g\mapsto \rho(g)S=(\pi(g)S^*)^*$ shows that $S\dd\pi(x_1)\cdots \dd\pi(x_k)$ is continuous in $(x_1,\cdots,x_k)$. Therefore \eqref{metrthmeq1} entails  that $\|\dd\pi(x)^n S \dd\pi(y)^m\|$ is bounded for $(x,y)$ in a neighborhood of $(0,0)\in \g^2$. Since $U(\g)$ is spanned by the elements of the form $x^k$, $x \in \g$, $k \in \N_0$, and 
$$U(\g) \times U(\g) \to B(\cH), \quad (D_1,D_2) \mapsto \dd\pi(D_1)S\dd\pi(D_2)$$
is bilinear, polarization implies that the $(n+m)$-linear map $f_{n,m}:\g^{n+m}\to B(\cH)$
$$ f_{n,m}(x_1,\dots,x_n, y_1,\dots,y_m) := \dd\pi(x_1)\cdots \dd\pi(x_n)S\dd\pi(y_1)\cdots \dd\pi(y_m)$$
is bounded near $0$ and therefore continuous for every $n,m\in \N$.

Next we show that (the unique extension to $\cH$ of) $\dd\pi(D_1)S\dd\pi(D_2)$ lies in $B(\cH)_c(\alpha)$ for every $D_1,D_2\in U(\g)$. The proof is inductive, namely, we assume that $\oline A\in B(\cH)_c(\alpha)$
where
$A:=\dd\pi(D_1)S\dd\pi(D_2)$, and we show that  for all $x,y\in\g$, the unique extensions to $\cH$ of $\dd\pi(x)A$ and $A\dd\pi(y)$ are in $B(\cH)_c(\alpha)$.
 Remark \ref{rmk-hadi-xy} and \cite[Lem.~2.9]{NSZ15} imply that 
$\oline A\in\cD(\oline{\dd\lambda}(x))$ for any $x\in\g$, and $\oline{\dd\lambda}(x)\oline A=\oline{\dd\pi}(x)\oline A$. Now an argument similar to Remark \ref{rem:1}(b) yields 
$\oline{\dd\pi}(x)\oline A \in B(\cH)_c(\alpha)$. 
Furthermore, $(\oline A)^*=A^*$ and it is straightforward to verify that 
\[A^*\big|_{\cH^\infty}=
\dd\pi(D_2^\dagger)S^*\dd\pi(D_1^\dagger),
\] where $\dagger$ is the principal anti-involution of $U(\g)$ defined by $x^\dagger:=-x$ for $x\in\g$. 
Since obviously $S^*$ is Schwartz, 
the operator $A^*$ is the unique extension to $\cH$ of the bounded operator $\dd\pi(D_2^\dagger)S^*\dd\pi(D_1^\dagger)$, hence for any $y\in\g$ we have  
 by Remark \ref{rmk-hadi-xy} that 
$A^*(\cH)\subseteq\cD(\oline{\dd\pi}(y)^2)$. Now \cite[Lem~2.8(a)]{NSZ15} yields boundedness of 
$\oline A\oline{\dd\pi}(y)^2$, and \cite[Lem~2.10]{NSZ15} implies that $\oline A\oline{\dd\pi}(y)\in B(\cH)_c(\alpha)$.

Next we observe that for $x,y\in \g$, the partial derivatives of
$$\R^2 \to B(\cH), \quad (t,s) \mapsto \pi(\exp(tx))\dd\pi(D_1)S\dd\pi(D_2)\pi(\exp(-sy))$$
exist and are continuous
(see \cite[Lemmas 2.9/10]{NSZ15} and Remark \ref{rmk-hadi-xy}, and recall from above that 
for 
$A=\dd\pi(D_1)S\dd\pi(D_2)$, the operator $\oline A\oline{\dd\pi}(y)^2$ is bounded). This yields $\dd\pi(D_1)S\dd\pi(D_2)\in \cD^1(\alpha)$ and 
$$\overline{\dd\alpha}(x,y)(\dd\pi(D_1)S\dd\pi(D_2))=\dd\pi(x)\dd\pi(D_1)S\dd\pi(D_2)-\dd\pi(D_1)S\dd\pi(D_2)\dd\pi(y).$$
Hence we can prove $S\in \cD^\infty(\alpha)$ by induction. The continuity of the maps $f_{n,m}$ and \cite[Cor.~1.7(ii)]{NSZ15} now implies that $S$ is a smooth vector for $\alpha$.
\end{prf}

Recall that $\mathcal S(\pi,\cH)$ denotes the space of Schwartz operators of a unitary representation
$(\pi,\cH)$.
The next proposition is an application of  
Theorem \ref{thm:2.4}.

\begin{prop}
\label{prp:funccalc}
Let $(\pi,\cH)$ be a smooth unitary representation of a Fr\'echet--Lie group $G$. Let $T\in \mathcal S(\pi,\cH)$. Assume that $T$ is a non-negative self-adjoint operator. Then $\sqrt{T}\in \mathcal S(\pi,\cH)$.  
\end{prop}

\begin{prf}
Since $\sqrt{T}$ is self-adjoint, by Theorem \ref{thm:2.4} it is enough to show that it is smoothing. Next choose $v\in\cH^\infty$ such that $\|v\|=1$. Let $\dagger$ denote the 
principal anti-involution of $U(\g)$, defined by $x^\dagger:=-x$ for $x\in\g$. Then
\begin{align*}
\|\sqrt{T}\dd\pi(D)v\|^2&=\langle \dd\pi(D)v,T\dd\pi(D)v
\rangle\\
&=\langle v, \dd\pi(D^\dagger)T\dd\pi(D)v
\rangle\leq \|\dd\pi(D^\dagger)T\dd\pi(D)\|.
\end{align*}
Thus the operator $\sqrt{T}\dd\pi(D)$ is bounded on $\cH^\infty$. From
 \cite[Thm 2.11]{NSZ15} it follows that $\sqrt{T}$ is smoothing.
\end{prf}

\section{Relation to literature on Schwartz operators}

Schwartz operators have also been studied 
in \cite{Pedersen} for nilpotent Lie groups, and more generally in \cite{Beltita}. Note that
from \cite[Thm.~2.11]{NSZ15} it follows that there is redundancy in the definitions given in \cite[Sec.~1.2]{Pedersen} and \cite[Def. 3.1]{Beltita}. 
From \cite[Thm.~3.1]{Beltita} it follows that 
smooth vectors of the $G\times G$-action on $B_2(\cH)$ are Schwartz operators. This is weaker than Theorem \ref{thm:2.4} above. 
Furthermore,   
\cite[Thm.~3.1]{Beltita} gets close to 
Proposition \ref{prop:1.6} and 
Corollary \ref{cor:1.8}(iii),
but in \cite{Beltita} it is not proved that being trace class is equivalent to nuclearity of the space of smooth vectors (see Proposition \ref{prop:3.18}). Finally, 
Proposition \ref{prop:1.6} implies that
what is proved in \cite[Cor.~3.1]{Beltita} for irreducible unitary representations of nilpotent Lie groups indeed holds for all trace class representations of general finite dimensional Lie groups.

\subsection*{The Schr\"{o}dinger representation}

In this section we investigate the connection between our results and those of \cite{KKW15} more closely. In particular, we will show that several of the  results of \cite{KKW15} are special cases of the results of our paper, when applied to the Schr\"{o}dinger representation. 

Let $(V,\omega)$ be a $2n$-dimensional real symplectic space and let
$H_{V,\omega}$ denote the  Heisenberg group associated to $(V,\omega)$, that is, $H_{V,\omega}:=V\times\R $ with the multiplication
\[
\textstyle(v,s)(w,t):=\left(v+w,s+t+\frac{1}{2}\omega(v,w)\right).
\]
Let
$\mathfrak h_{V,\omega}$ denote 
the Lie algebra  of $H_{V,\omega}$, and 
let $U(\mathfrak h_{V,\omega})$ denote the universal enveloping algebra of $\mathfrak h_{V,\omega}$. 
By the Stone--von Neumann Theorem, to every nontrivial 
unitary character $\chi:\R\to\C^\times$  we can associate  a unique irreducible unitary representation $\pi_\chi$ of $H_{V,\omega}$ for which the center acts by $\chi$. In the Schr\"{o}dinger realization, $\pi_\chi$ acts on the Hilbert space
$\cH:=L^2(Y,\mu)$, where $V=X\oplus Y$ is a polarization of $V$, and $\mu$ is the Lebesgue measure on $Y\cong \R^n$. The action of $\pi_\chi$ is given by
\[
(\pi_\chi(x,0)\phi)(y)
:=\chi(\omega(x,y))\phi(y),\quad 
(\pi_\chi(y_0,0)\phi)(y):=\phi(y-y_0),
\]
and $\pi_\chi(0,t)\phi:=\chi(t)\phi$,
where $x\in X$, $\phi\in L^2(Y,\mu)$, $y,y_0\in Y$, and $t\in\R$. 
The following result is a special case of the general theory of unitary representations of nilpotent Lie groups
(e.g., see 
\cite{Howe80}).

\begin{prop}
\mlabel{prop:pichi}
The representation $\pi_\chi$ is  trace class, the space of smooth vectors of $\pi_\chi$ is the Schwartz space $\cS(Y)$, and $\dd\pi(U(\mathfrak h_{V,\omega}))$ is equal to the algebra of polynomial coefficient differential operators on $Y$.
\end{prop}
From Proposition \ref{prop:pichi} it follows that the operators defined in 
\cite[Def.~3.1]{KKW15} are  the Schwartz operators for $\pi_\chi$ in the sense of our paper. From 
Corollary~\ref{cor:1.8} it follows that $\mathcal S(\pi_\chi,\cH)$ is the space of smooth vectors for the action of 
$H_{V,\omega}\times H_{V,\omega}$ on $B_2(\cH)$, and therefore it can be equipped with a canonical  Fr\'echet topology. It is straightforward to verify that this Fr\'echet topology is identical to the one described in \cite[Prop. 3.3]{KKW15}. 

Next we show  that if $S,T\in \mathcal S(\pi_\chi,\cH)$, then 
$SAT\in \mathcal S(\pi_\chi,\cH)$ for every $A\in B(\cH)$. This is proved in 
\cite[Lemma~3.5(b)]{KKW15}, but the argument that will be given below applies to any trace class representation. From
Proposition \ref{prop:1.7} it follows that 
the map $G\to B_2(\cH)$ given by $g\mapsto \lambda(g)S$ is smooth. Since the bilinear map
\[
B_2(\cH)\times B(\cH)\to B_2(\cH)\ ,\ 
(P,Q)\mapsto PQ
\] is continuous, the map $g\mapsto \lambda(g)SAT$ is also smooth. Thus by Proposition \ref{prop:1.7},
the operator $SAT$ is smoothing. A similar argument shows that $T^*A^*S^*$ is also smoothing, and Corollary \ref{cor:1.8}
implies that $SAT$ is Schwartz. 

Proposition \ref{prop:1.6} implies  
that every Schwartz operator for $\pi_\chi$ is trace class. This is also obtained in \cite[Lemma~3.6]{KKW15}. Theorem \ref{thm:2.4} applied to $\pi_\chi$ gives \cite[Thm.~3.12]{KKW15}.
Proposition \ref{prp:funccalc} implies \cite[Prop. 3.15]{KKW15}, and Proposition \ref{prp:linfunc} implies 
\cite[Prop. 5.12]{KKW15}.

It is possible that the relation between the Weyl transform and Schwartz operators that is investigated in \cite[Sec.~3.6]{KKW15} is a special case of  more general results in the spirit of our paper, at least for nilpotent Lie groups. Finally, it is worth mentioning that the paper 
\cite{Be11} studies (among several other things) 
the class of representations of infinite dimensional Lie groups with the property that their space of 
smooth vectors is nuclear. By Proposition \ref{prop:3.18}, when $G$ is finite dimensional this condition is equivalent to the representation being trace class. In the infinite dimensional case this is an interesting class of representations which deserves further investigation. 
  We hope to come back to these problems in the near future.

\end{document}